\documentclass{article}
\usepackage{texdraw}
\newtheorem{thm}{Theorem}[section]

\newtheorem{prop}[thm]{Proposition}

\newcommand{\be}{\begin{equation}}
\newcommand{\ee}{\end{equation}}
\newcommand{\ben}{\begin{enumerate}}
\newcommand{\een}{\end{enumerate}}
\newcommand{\beq}{\begin{eqnarray}}
\newcommand{\eeq}{\end{eqnarray}}
\newcommand{\beqn}{\begin{eqnarray*}}
\newcommand{\eeqn}{\end{eqnarray*}}

\newcommand{\pa}{\partial}

\newcommand{\pxi}{ {\pa \over \pa x^i}}

\newcommand{\qed}{\hspace*{\fill}Q.E.D.} %Use at end of proof

\title{A Family of Einstein Randers Metrics}

\author{B. Najafi and A. Tayebi}
\begin{document}

\maketitle

\begin{abstract}
Using Hawking Taub-NUT metric  $g_a$ on $\mathbf{R}^4$, where $a$ is a positive real number and finding a 4-parameter family
of Killing vector fields of $(\mathbf{R}^4,g_a)$, we construct a 5-parameter family of Einstein Randers metrics with non-constant
flag curvature.
\end{abstract}
\section{Introduction}
The Taub-Newman-Unti-Tamburino (Taub-NUT) metrics were found by Taub
\cite{T} and extended by Newman-Unti-Tamburino \cite{NTU}. The Euclidean Taub-NUT
metric has lately attracted much attention in physics. Hawking has suggested
that the Euclidean Taub-NUT metric might give rise to the gravitational analog
of the Yang-Mills instanton \cite{Ha}. This metric is the space part of the line element of
the celebrated Kaluza-Klein monopole of Gross and Perry and Sorkin.  The Taub-NUT
family of metrics is also involved in many other modern studies in physics like
strings, membranes, etc.

On the other hand, a rich and in some sense one of the best classes of Riemannian metrics is Einstein metrics, i.e, the solutions of the equation $Ric_{ij}-\lambda g_{ij}=0$ for some scalar $\lambda$. Einstein manifolds are not only interesting in physics but are also related to many important topics of Riemannian geometry such as self-dual manifolds of dimension four. In dimension three, being Riemannian Einstein metric is equivalent to being of constant sectional curvature. S. S. Chern conjectured that this equivalency is not true for Finslerian Einstein metrics. In \cite{BRo1}, D. Bao and C. Robles classified  Einstein Randers metrics in every dimension and showed that in dimension three, any Einstein Randers metric is of constant flag curvature (a natural generalization of sectional curvature) and vice versa. But for general Finslerian Einstein metrics, Chern's conjecture remains unsolved.

After classifying 3-dimensional Einstein metrics of Randers type, it is natural to classify 4-dimensional Einstein Randers metrics. The reasonable step is constructing as many as possible concrete examples. In this paper, we construct a 5-parameter family of  4-dimensional Randers Einstein metrics with non-constant flag curvature. This family contains the family of Randers Einstein metrics constructed in \cite{GMZ} as special case.
\section{Preliminaries}
Let $F$ be a Finsler metric on an $n$-dimensional manifold $M$, and $TM_0$ be its slant tangent space. In a local coordinate $(x^i,y^i)$, the vector filed $G=y^i\frac{\pa}{\pa x^i}-2G^i\frac{\pa}{\pa y^i}$ is a global vector field on $TM_0$,
where $G^i=G^i(x,y)$ are local functions on $TM_0$ satisfying $G^i(x,\lambda y)=\lambda^2 G^i(x,y)$,  for all $\lambda>0$. The vector field $G$ is called the associated spray to $(M,F)$.

The  Riemann tensor can be viewed as a family of endomorphisms on tangent spaces: ${\bf R}_y (u): = R^i_{\ k}(x,y) u^k \pxi|_x,$
where $u=u^i\pxi|_x \in T_xM$. The coefficients $R^i_{\ k}=R^i_{\ k}(x,y)$ are given by
\[
R^i_{\ k} = 2 {\pa G^i\over \pa x^k}-y^j{\pa^2 G^i\over \pa x^j\pa y^k}
+2G^j {\pa^2 G^i \over \pa y^j \pa y^k} - {\pa G^i \over \pa y^j}{\pa G^j \over \pa y^k}.  \label{Riemann}
\]
For a tangent plane $P=span(y,u)\subset T_pM$, the \textit{flag curvature} of the flag
$P$ with the flag pole $y$ is defined by
\[
K(P,y)={{g_y({\bf R}_y (u),u)} \over {g_y(y,y)g_y(u,u)-g_y(u,y)g_y(u,y)}} .
\]
If $F$ is Riemannian, i.e., $g=g_y$ is independent of $y$, then the
flag curvature $K(P,y)=K(P)$ is independent of $y\in P$ and it is
the sectional curvature of $P$. Thus the flag curvature is an
analogue of the sectional curvature in Riemannian geometry. F is
said to be of \textit{constant flag curvature} $\lambda$, if for any non-zero tangent vector $y \in
T_xM$ and any flag $P$ containing $y$,
$K(P,y)=\lambda$ is constant, or equivalently,
\[
{\bf R}_y =\lambda F^2(y)\{I-g_y(y,.)y\}, \quad y \in T_xM, \quad x \in M
\]
where $I: T_xM \rightarrow T_xM$  denotes the identity map and
$g_y(y,.)={{1}\over{2}}[F^2]_{y^i}dx^i$.

The trace of
${\bf R}_y$ is denoted by ${\bf Ric}(y)$ and called the {\it Ricci curvature}.
Hence, $Ric(x,y)=R^i_i(x,y)$. A Finsler metric $F$ is said to be an \textit{Einstein metric} if there is a scalar function
$K=K(x)$ on $M$ such that $Ric=(n-1)KF^2$.

In \cite{W}, Wood constructs the following 4-dimensional Riemannian metric. Let $(N^3,h)$ be an oriented 3-dimensional Riemannian manifold with
constant curvature. Let $M^4:=\mathbf{R}\times N^3$ and $\phi:M^4\rightarrow N^3$ be the projection on the second factor. Suppose that $u$ is a positive smooth function and $A$ a 1-form on $N^3$, then the tensor $g$ on
$M^4$ given by
\begin{equation} \label{harmonic metric}
g:=u\phi^*(h)+ u^{-1}(dt+A)^2
\end{equation}
is a Riemannian metric on $M^4$.
In \cite{GHa}, Gibbons and Hawking prove that $(M^4,g)$ is an Einstein manifold if and only if the monopole equation $du=*dA$ holds and
$(N^3,h)$ is flat, in which case $g$ is Ricci-flat.

Let us consider the following harmonic function on $\mathbf{R}^3-\{0\}$
\[
u_a(x):=\frac{1}{4}\big(\frac{1}{|x|}+a\big),
\]
where $a$ is a non-negative real number. Then the above construction gives the Hawking Taub-NUT metric $g_a$ $(a>0)$ or
the standard metric $g_0$ ($a=0$). A direct computation yields the following explicit formula for $g_a$
\begin{equation}
g_a=(a|x|^2+1)g_0-\frac{a(a|x|^2+2)}{a|x|^2+1}\omega^2,
\end{equation}
where $g_0$ is the standard metric and
\begin{equation}\label{form}
\omega=-x^2dx^1+x^1dx^2-x^4dx^3+x^3dx^4
\end{equation}
is a 1-form on $\mathbf{R}^4$ (see \cite{W} for more details). Hence the metric $g_a$ extends to the whole $\mathbf{R}^4$. \\

A Randers metric on a manifold $M$ is a Finsler metric in the following form
\begin{equation} \label{randers}
F=\alpha+\beta,
\end{equation}
where $\alpha=\sqrt{a_{ij}(x)y^iy^j}$ is a Riemannian metric and $\beta=b_i(x)y^i$ a 1-form on $M$. A nice description of Randers metrics is given by Zermelo navigation representation. Given a Riemannian metric $g=\sqrt{g_{ij}y^iy^j}$ and a vector field $V=V^i\frac{\pa}{\pa x^i}$, set
$V_i:=g_{ij}V^i$ and $||V||_g:=\sqrt{V_iV^i}$. Then the pair $(g,V)$ with $||V||_g<1$ is called the \textit{Zermelo navigation representation} (or simply navigation representation) of the Randers metric (\ref{randers}) if and only if the following relations hold
\begin{equation}\label{zermelo 1}
a_{ij}=\frac{g_{ij}}{1-||V||^2_g}+\frac{V_iV_j}{(1-||V||^2_g)^2}, \quad\quad b_i=-\frac{V_i}{1-||V||^2_g}.
\end{equation}\\
It is easy to see that, if $(g,V)$ is the navigation representation of $F=\alpha+\beta$, then the following hold
\begin{equation}
F(x,y)=g(x,y-F(x,y)V(x)).
\end{equation}
Many curvature properties of Randers metrics can be expressed in a geometrical way via their navigation representations. D. Bao and C. Robles
 use this method and characterize Einstein Randers metrics as follows \cite{BR}.
\begin{prop} \label{robles 1} \emph{Suppose that a Randers metric $F=\alpha+\beta$ has the navigation representation $(g,V)$ on a manifold $M$. Then $F$ is an Einstein metric with $Ric(x,y)=(n-1)K(x)F^2(x,y)$ if and only if
\begin{description}
  \item[i)] $V$ is a homothetic vector field of $g$ with $\mathcal{L}_V(g)=-4cg$, where $\mathcal{L}_V$ denotes the Lie derivative along $V$.
  \item[ii)] $g$ is Einstein metric with $\widetilde{Ric}=(n-1)\{K(x)+c^2\}g$, where $\widetilde{Ric}$ is
the Ricci curvature of $g$. In particular, $K(x)=const$ if $n\geq 3$.
\end{description}}
\end{prop}

\section{Construction of the Einstein Finsler Metrics}

Rewrite the Hawking Taub-NUT metric $g_a$ as
\[
g_a=g_{ij}dx^idx^j,
\]
where
\[
(g_{ij})=\left(
           \begin{array}{cccc}
             B-A(x^2)^2 & Ax^1x^2 & -Ax^2x^4 &Ax^2x^3 \\
            Ax^1x^2 & B-A(x^1)^2 & Ax^1x^4 & -Ax^1x^3 \\
             -Ax^2x^4 & Ax^1x^4 & B-A(x^4)^2& Ax^3x^4 \\
             Ax^2x^3 & -Ax^1x^3 & Ax^3x^4 & B-A(x^3)^2\\
           \end{array}
         \right)
\]
and $B=B(x)=a|x|^2+1$, $\qquad A=A(x)=a\Big(1+\frac{1}{B}\Big)$.

The key observation is this: Suppose that $f:\mathbf{R}^4\rightarrow \mathbf{R}^4$ is an isometry of $g_0$, preserves the 1-form $\omega$ given by
(\ref{form}) and $f(0)=0$. Then, it is easy
to see that $f$ is also an isometry of $g_a$. Now we classify those Killing vector fields of $g_0$ which are also
Killing vector fields with respect to $g_a$ and have the origin as a critical point. First, recall the following characterization theorem \cite{BRS}.
\begin{prop}\label{prop 1}\emph{  Let $X=X^i\pxi$ be a vector field on $\mathbf{R}^4$. Then $X$ is a Killing vector field of $g_0$
if and only if $X^i=Q^i_jx^j+C^i$, where  $(Q^i_j)$ is an anti-symmetric matrix,  $(C^i)$ is a constant vector.}
\end{prop}

Using Proposition \ref{prop 1}, we find a 4-parameter family of Killing vector fields of $g_a$, and consequently a 5-parameter family of Einstein Randers metric with non-constant flag curvature. More precisely, with the notation of Proposition \ref{prop 1}, we have the following.
\begin{thm}\label{thm 1} Let $X=X^i\pxi$ be a Killing vector field on $(\mathbf{R}^4,g_0)$. Then $X$ is a Killing vector field on $(\mathbf{R}^4,g_a)$ if and only if $(C^i)=(0)$ and the matrix $(Q^i_j)$ is in the following form
\[
\left(
  \begin{array}{cccc}
    0 & m & s & r \\
    -m & 0 & -r & s \\
    -s & r & 0 & n \\
    -r & -s & -n & 0 \\
  \end{array}
\right),
\]
where $m$, $n$, $r$ and $s$ are constant real numbers.
\end{thm}
{\bf{Proof.}}
Using the geometric meaning of Lie derivation of vector fields, the one parameter groups of the vector field $X=X^i\pxi$ on $\mathbf{R}^4$ preserves
the 1-form $\omega$ if and only if $\mathcal{L}_X(\omega)=0$. A direct computation implies that $\mathcal{L}_X(\omega)=0$ is equivalent to the following system of PDEs
\begin{equation}
\left\{
  \begin{array}{ll}
    -X^2-x^2 \frac{\pa X^1}{\pa x^1}+x^1\frac{\pa X^2}{\pa x^1}-x^4\frac{\pa X^3}{\pa x^1}+x^3\frac{\pa X^4}{\pa x^1}& \hbox{=0,} \\
    -x^2\frac{\pa X^1}{\pa x^2}+X^1+x^1\frac{\pa X^2}{\pa x^2}-x^4\frac{\pa X^3}{\pa x^2}+x^3\frac{\pa X^4}{\pa x^2} & \hbox{=0,} \\
    -x^2\frac{\pa X^1}{\pa x^3}+x^1\frac{\pa X^2}{\pa x^3}-X^4-x^4\frac{\pa X^3}{\pa x^3}+x^3\frac{\pa X^4}{\pa x^3} & \hbox{=0,} \\
    -x^2\frac{\pa X^1}{\pa x^4}+x^1\frac{\pa X^2}{\pa x^4}-x^4\frac{\pa X^3}{\pa x^4}+X^3+x^3\frac{\pa X^4}{\pa x^4}& \hbox{=0.}
  \end{array}
\right.
\end{equation}
By $\frac{\pa X^i}{\pa x^i}=\sigma$ and  $\frac{\pa X^i}{\pa x^j}=Q^i_j$ for distinct $i$ and $j$, we get the result. \qed
\newpage
Using Proposition \ref{robles 1}, we complete the construction in three cases.
\subsection{r=0, s=0.} If we put $Q^1_2=-m$ and $Q^3_4=-n$, where $m$ and $n$ are constant real numbers, then the vector field $X$ reduces to $W_{m,n}$ which is introduced by E. Guo, X. Mo, and X. Zhang in \cite{GMZ}. Using navigation representation $(g_a,W_{m,n})$, a 3-parameter family of  Einstein Randers metric with non-constant flag curvature is constructed on the $\Omega=\{x\in\mathbf{R}^4  \, |\, f(x)<1\}$, where
\[
f(x):=|W_{m,n}|_{g_a}=\frac{|x|^2}{1+a|x|^2}(p+2a|m-n||x|^2+a^2|m-n||x|^4).
\]
For more details see \cite{GMZ}.\\

\subsection{m=0, n=0, s=0.}
In this case, the vector field $X$ reduces to the vector field $V_r$ given by
\begin{equation}\label{new vect}
V_{r}:=r(x^4\pa _1-x^3\pa_2+x^2\pa_3-x^1\pa_4),
\end{equation}
where $r$ is a constant real number. Here, we construct Einstein Randers metric with non-constant flag curvature by using navigation representation $(g_a,V_r)$. First, we find the sufficient condition producing Einstein-Finsler metrics, i.e., $||V_r||_{g_a}<1$.\\

We have $V^1=rx^4$, $V^2=-rx^3$, $V^3=rx^2$, and $V^4=-rx^1$. A direct computation yields
\[
V_i:=g_{ij}V^j=\mu V^i,
\]
where $\mu=\frac{1}{a|x|^2+1}$. Therefore, we get
\begin{equation} \label{norm of V}
||V_r||^2_{g_a}=V_iV^i=\mu r^2|x|^2.
\end{equation}

If $r^2>a$, then define $\Omega_{a,r}=\{x\in\mathbf{R}^4  \, |\,|x|<\frac{1}{\sqrt{r^2-a}} \}$. Otherwise, define $\Omega_{a,r}=\mathbf{R}^4$.

\begin{thm}
Let $F=\alpha+\beta$ be a Randers metric be given by the navigation representation $(g_a,V_r)$. Then $F$ is an Einstein metric with non-constant flag curvature on $\Omega_{a,r}$.
\end{thm}
{\bf{Proof.}} From (\ref{norm of V}), we get  $||V_r||_{g_a}<1$ on $\Omega_{a,r}$. Hence $F$ is actually a Randers metric on $\Omega_{a,r}$.

In \cite{BRo2}, D. Bao and C. Robles prove that a Randers metric with navigation representation $(g,V)$ is of constant flag curvature if and only if
$g$ is of constant sectional curvature and $V$ is a homothetic vector field with respect to $g$. The Taub-NUT metric is a Ricci flat metric on $\mathbf{R}^4$ which is not flat \cite{Ha}. Hence $g_a$ is not of constant sectional curvature. Therefore $F$ is not of constant flag curvature.

It is sufficient to prove that $F$ is Einstein metric.  Note that the Hawking Taub-NUT metric
$g_a$ for all $a>0$ is an Einstein metric \cite{W}. From Theorem \ref{thm 1}, we see that $V_r$ is a Killing vector field with respect to
$g_a$. As an immediate consequence of Proposition \ref{robles 1}, we conclude that $F$ is an Einstein metric. This completes the proof. \qed

\subsection{m=0, n=0, r=0.}
In this case, the vector field $X$ reduces to the vector field $U_s$ given by
\begin{equation}\label{new vect}
U_{s}:=s(x^3\pa _1+x^4\pa_2-x^1\pa_3-x^2\pa_4),
\end{equation}
where $s$ is a constant real number. We have
\[
U^1=sx^3, U^2=sx^4, U^3=-sx^1,  U^4=-sx^2.
\]
Using a Maple program, we get
\begin{equation} \label{norm of U}
||U_s||^2_{g_a}=\mu s^2\{(1+8a^2 \prod_{i=1}^{4} x^i)|x|^2+2a|x|_4^{\,4}+a^2|x|_6^{\,6}+16a\prod_{i=1}^{4} x^i+f(x)\},
\end{equation}
where $|.|_k$ denotes the $k$-norm on $\mathbf{R}^4$ and the scalar function $f$ is given by
\begin{eqnarray}
            f(x)=&-& a^2(x_1^2x_4^4+x_2^4x_3^2+x_2^2x_3^4+x_1^4x_4^2) \nonumber\\
            &+& 3a^2(x_1^2x_2^4+x_4^2x_2^4+x_1^2x_3^4+x_3^2x_4^4+x_3^4x_4^2+x_1^4x_2^2+x_2^2x_4^4+x_3^2x_1^4) \nonumber\\
            &+& 4a(x_1^2x_3^2-x_2^2x_3^2+x_2^2x_1^2+x_4^2x_3^2-x_4^2x_1^2+x_4^2x_2^2) \nonumber\\
            &+& 2a^2(x_1^2x_3^2x_4^2+x_2^2x_3^2x_4^2+x_1^2x_2^2x_3^2+x_1^2x_2^2x_4^2)\nonumber,
\end{eqnarray}
in which $x_i^k$ means $(x^i)^k$.\\\\
Set
$$\widetilde{\Omega}_{a,s}:=\{x\in \mathbf{R}^4|\,\, ||U_s||_{g_a}<1 \}.$$
Then we have the following theorem.
\begin{thm}
Let $F=\alpha+\beta$ be a Randers metric given by the navigation representation $(g_a,U_s)$. Then $F$ is an Einstein metric with non-constant flag curvature on $\widetilde{\Omega}_{a,s}$.
\end{thm}

These concrete examples of four dimensional Randers Einstein metrics show that this class of Randers metrics is large and of course far from being
completely classified.

\bigskip
Behzad Najafi\\
Faculty of Science, Department of Mathematics\\
Shahed University of Tehran\\
Tehran. Iran\\
Email: najafi@shahed.ac.ir
\bigskip
  \\[1.5ex]
Akbar Tayebi\\
Faculty of Science, Department of Mathematics\\
Qom University\\
Qom. Iran\\
Email:\  akbar.tayebi@gmail.com

\end{document}